\documentclass[runningheads,envcountsame]{llncs}

\usepackage[T1]{fontenc}
\usepackage[utf8]{inputenc}
\usepackage[english]{babel}
\usepackage{amsmath,amssymb,mathtools}
\usepackage{microtype}
\usepackage[hidelinks]{hyperref}

\newcommand{\Lang}{\mathcal L}
\newcommand{\LangQ}{\mathcal L_{\mathbb Q}}
\newcommand{\Mod}{\operatorname{Mod}}

\newcommand{\K}{\mathcal K}
\newcommand{\val}[2]{\lVert #1\rVert_{#2}}
\newcommand{\forcesval}[2]{[#1]^{#2}}

\newcommand{\preck}[1]{\preccurlyeq_{#1}}

\begin{document}

\title{Chang's Omitting Types Theorem in {\L}ukasiewicz Logic}
\titlerunning{Chang's Omitting Types Theorem in {\L}ukasiewicz Logic}

\author{Luca Spada}
\authorrunning{L. Spada}
\institute{Department of Mathematics, University of Salerno,\\ Piazza Renato Caccioppoli, 2, 84084 Fisciano (SA), Italy\\
\email{lspada@unisa.it}}

\maketitle

\begin{abstract}
We prove Chang-style omitting types theorems for standard first-order
{\L}ukasiewicz logic.  Starting from a prescribed model $M$, we obtain an extension that omits higher-complexity types while preserving lower-level
types already realized in $M$. If $T$ is a $\Pi_{n+2}$-theory, every model of $T$ has an \emph{exact} extension that preserves its realized $\Sigma_{n+1}$-types and omits every $\Sigma_{n+2}$-type not generated by a $\Sigma_{n+1}$-type.  An iterated construction also omits infinity-universal types, whose unbounded quantifier complexity cannot eventually be generated by the corresponding lower-complexity fragments of a single type.
The main tool is a refinement of Robinson's infinite forcing sensitive to the quantifier complexity of formulas.

This is an extended version of \cite{Spada26}.

\end{abstract}

\section{Introduction}

The classical omitting types theorem grew out of the diagram method developed
by Henkin and Robinson.  Its foundational form is implicit in
Henkin's generalization of $\omega$-consistency \cite{Henkin1954} and in Orey's subsequent work \cite{Orey1956}; see also Keisler's historical discussion and forcing formulation \cite{Keisler1973}. In its familiar countable form, a countable theory together with a
countable family of nonprincipal types has a countable model omitting all of
them.  The result became a basic tool for constructing models with prescribed
failures of realization.  Robinson's infinite forcing and the general forcing
framework developed later by Keisler made explicit the close relation
between diagram constructions, generic models, and type omission
\cite{Keisler1973}.

Chang's theorem on prenex formulas \cite{Chang1967} represents a different branch of this
history.  Rather than merely producing some model of a theory that omits a prescribed type, it begins with a given model $M$ and constructs an extension in which types are realized or omitted depending on their quantifier complexity: for a $\Pi_{n+2}$-theory, the
extension preserves the $\Sigma_{n+1}$-types already realized in $M$ while
omitting the $\Sigma_{n+2}$-types that are not generated by any
$\Sigma_{n+1}$-type. Lablanquie later recovered Chang's result from a version of Robinson's forcing accounting for the quantifier complexity of formulas \cite{Lablanquie1978}.  

In the context of real-valued model theory, omitting types theorems have already been obtained.  In \cite[Theorem~3.6]{CaicedoIovino2014}, the
$\lambda$-omitting types property for continuous logic is proved for every infinite cardinal $\lambda$; suitable restrictions of that framework include
{\L}ukasiewicz-Pavelka logic.  Independently, in \cite[Theorem~4.9]{Eagle2014} it is proved that, for countable fragments of an infinitary $[0,1]$-valued logic, countably many nonprincipal types can be omitted.  Both
approaches are topological and use Baire-category arguments.  They accommodate broad families of nonprincipal types and, in the setting of \cite{CaicedoIovino2014}, arbitrary infinite cardinalities.  In spirit, they are closer to the omitting types theorems of Henkin and Orey.

To obtain a result similar to Chang's in the context of {\L}ukasiewicz logic, we follow the approach of Keisler and Lablanquie based on infinite forcing. Infinite forcing was extended to standard first-order {\L}ukasiewicz logic in \cite{DiNolaGeorgescuSpada2008}.  Model-theoretic forcing has also been developed for metric structures by Ben Yaacov and Iovino in \cite{BenYaacovIovino2009}.
Unlike the above Baire-category arguments, the proof given here is
\emph{constructs} a model starting from a
prescribed one. This should not be understood as an
effective procedure. The proof still uses standard satisfiability compactness, choices of suitable extensions, and transfinite recursion.

There is a specific obstruction to applying the classical Lablanquie's proof.  In a classical structure, knowing which formulas are true determines the truth values of all formulas.  In a {\L}ukasiewicz structure it does not: two formulas may both fail to have value one while having different values in $[0,1]$.
We repair the argument in two steps.  First, we replace truth-one
substructures at level $n$ by \emph{degree-preserving} ones.  Rational
truth constants provide a syntactic diagram for this semantic relation.  The
use of the rational Pavelka expansion is conservative over  {\L}ukasiewicz logic \cite{HajekParisShepherdson2000}: every standard structure in the original language has a unique canonical expansion, so reducts return us to the original setting.
Second, the classical Lablanquie lemma must be adapted to the many-valued setting.  Unsatisfiability yields, by compactness, a uniform gap below truth degree one; a suitable finite power of a formula turns that gap into zero and therefore produces a positive lower bound rather than an exact truth-one consequence.  For every rational $r\in[0,1)$ we obtain a lower-level formula forcing the original formula to have value at least $r$.  The collection of all these bounds forces value one and is sufficient for omitting types.

These modifications yield two relative omitting types theorems.  The first is
the direct analog of Chang's prenex theorem; the second iterates generic
extensions through the hierarchy and omits infinity-universal
types.  The latter notion is formulated semantically, in direct analogy with Lablanquie's classical condition.

The paper is organized as follows.  Section~\ref{sec:preliminaries} contains preliminaries fixing semantics, introducing degree diagrams, and proving the chain and extension lemmas.  In Section~\ref{sec:forcing} we define $n$-forcing and construct generic extensions.  Section~\ref{sec:mainlemma} proves the approximate Lablanquie lemma. The omitting types theorems are derived in Section~\ref{sec:firstchang}.

\section{Standard semantics and the prenex hierarchy}
\label{sec:preliminaries}

We refer to
\cite{CintulaHajekNoguera2011,Hajek1998,HajekCintula2006}
for background on first-order many-valued logics, and to
\cite{CignoliOttavianoMundici2000} for MV-algebras and
{\L}ukasiewicz logic.

We work with nonempty structures for a relational first-order language
$\Lang$ with constants. A \emph{standard} $\Lang$-structure $M$ has a
nonempty domain, interprets each constant as an element of that domain,
and interprets every $k$-ary predicate symbol $P$ as a function $P^M\colon M^k\longrightarrow[0,1]$.
Thus, formulas take their truth values in the standard MV-algebra
$[0,1]_{\text{\L}}$. In the more general algebraic semantics for predicate
{\L}ukasiewicz logic, truth values may instead belong to an arbitrary
MV-chain or MV-algebra. Throughout this paper, however, \emph{model},
\emph{satisfiable}, and $\models$ always refer to the standard
$[0,1]$-valued semantics unless explicitly stated otherwise.

Equality, when present, is interpreted \emph{crisply}: for $a,b\in M$,
\[
 \val{a=b}{M}=
 \begin{cases}
  1,&\text{if }a=b,\\
  0,&\text{if }a\neq b.
 \end{cases}
\]
If the original language is equality-free, it may safely be expanded by crisp equality. Function symbols can similarly be included by requiring substructures to be closed under their interpretations.

For $u,v\in[0,1]$, let $u\Rightarrow v=\min(1,1-u+v)$ denote the {\L}ukasiewicz residuum. Taking $\neg$, $\to$, and $\exists$ as primitive connective, truth values are defined recursively by
\begin{align*}
 \val{\neg\varphi[s]}{M}
   &=1-\val{\varphi[s]}{M},\\
 \val{(\varphi\to\psi)[s]}{M}
   &=\val{\varphi[s]}{M}\Rightarrow\val{\psi[s]}{M}\\
 \val{\exists x\,\varphi[s]}{M}
   &=\sup\{\val{\varphi[s(x/a)]}{M}\mid a\in M\}.
\end{align*}
We say that $M$ \emph{satisfies} $\varphi[s]$, and write
$M\models\varphi[s]$, when $\val{\varphi[s]}{M}=1$.
This is sometimes called \emph{truth-one satisfaction}. Notice that
$M\models\exists x\,\varphi(x)$ does not necessarily mean that some
$a\in M$ satisfies $\varphi(a)$ to degree one: the supremum may be one
without being attained. 
As usual, we introduce the derived connectives
\[
 \varphi\oplus\psi=\neg\varphi\to\psi,\quad
 \varphi\odot\psi=\neg(\neg\varphi\oplus\neg\psi),\quad
 \forall x\,\varphi=\neg\exists x\,\neg\varphi,\quad
 \varphi^m=
 \underbrace{\varphi\odot\cdots\odot\varphi}_{m\text{ times}}.
\]

Finally, the canonical rational expansion $\LangQ$ contains a truth
constant $\bar q$ for every $q\in\mathbb Q\cap[0,1]$, always interpreted
by $q$. Every standard $\Lang$-structure has a unique canonical
$\LangQ$-expansion. We write $\LangQ(M)$ for the expansion of $\LangQ$ by a constant for each element of $M$.
A theory is a set of sentences. 

Recall that satisfiability compactness holds for predicate {\L}ukasiewicz logic:
\begin{theorem}[{\cite{Belluce_Chang_1963}, see also \cite[Theorem 5.4.24]{Hajek1998}}]\label{thm:compactness}
A theory has a standard model whenever each of its finite subsets has one.
\end{theorem}

\begin{remark}
This theorem does not say that truth-one semantic consequence is compact. For example, let $p$ be an atomic sentence and let
\[
 \Delta=\left\{p\to\overline{1/m}\mid m\geq2\right\}.
\]
Since $\val{p\to\overline{1/m}}{N}=1$ exactly when
$\val{p}{N}\leq1/m$, every standard structure satisfying all formulas in $\Delta$ has $\val{p}{N}=0$. Hence $\Delta\models\neg p$. However, no finite subset of $\Delta$ has this consequence: finitely many of the displayed conditions impose only a positive upper bound on $\val{p}{N}$, which may therefore be chosen strictly positive. 
\end{remark}

The next elementary observation will be used repeatedly.
\begin{lemma}\label{lem:power}
For every standard structure $M$ and every sentence $\varphi$ of its
expanded language, if $\val{\varphi}{M}<1$ then there exists some $m\geq1$ such that $\val{\varphi^m}{M}=0$.
\end{lemma}

\begin{proof}
Put $a=\val{\varphi}{M}$. For $a\in[0,1]$,
$a^{m}=\max(0,ma-(m-1))$. Choose $m$ with $m(1-a)\geq1$.
\qed
\end{proof}

The following lemma can be used as a semantic substitute for the deduction theorem.
\begin{lemma} \label{lem:uniformseparation}
Let $S$ be a theory, and let $\alpha_1,\ldots,\alpha_k$ be finitely many
sentences.  If $S\cup\{\alpha_1,\ldots,\alpha_k\}$ is not satisfiable, then there is an integer $m\geq1$ such that no model $N\models S$ satisfies
\[
 \val{\alpha_i}{N}\geq1-\frac1m
 \qquad(1\leq i\leq k).
\]
\end{lemma}

\begin{proof}
Suppose otherwise.  Consider the theory
\[
 S\cup
 \{\bar q\to\alpha_i \mid q\in\mathbb Q\cap[0,1),\ 1\leq i\leq k\}.
\]
Each of its finite subsets is satisfiable. Indeed, choose $m$ so large that
$1-1/m$ is greater than every rational $q$ occurring in that finite subset.
By the absurdum hypothesis, there is a model $N\models S$ in which every
$\alpha_i$ has value at least $1-1/m$; this model satisfies the finite set.
Therefore, Theorem~\ref{thm:compactness} gives a model of the displayed
theory. In that model, the rational constants force every $\alpha_i$ to have
value one, contrary to the hypothesis.
\qed
\end{proof}

\begin{lemma} \label{lem:prenex}
If $x$ is not free in $\psi$, the following logical equivalences hold over all
standard {\L}ukasiewicz structures:
\begin{align*}
 \neg\exists x\,\varphi&\equiv\forall x\,\neg\varphi,&
 \neg\forall x\,\varphi&\equiv\exists x\,\neg\varphi,\\
 (\exists x\,\varphi)\to\psi&\equiv
      \forall x\,(\varphi\to\psi),&
 (\forall x\,\varphi)\to\psi&\equiv
      \exists x\,(\varphi\to\psi),\\
 \psi\to\exists x\,\varphi&\equiv
      \exists x\,(\psi\to\varphi),&
 \psi\to\forall x\,\varphi&\equiv
      \forall x\,(\psi\to\varphi).
\end{align*}
Consequently every formula is equivalent over standard semantics to a
prenex formula.
\end{lemma}

\begin{proof}
Let $S\subseteq[0,1]$ be the nonempty set of values of $\varphi(a)$ as $a$ ranges over the domain, and let $b$ be the value of $\psi$.  Continuity and the monotonicity or antitonicity of the {\L}ukasiewicz residuum give
\begin{align*}
 (\sup S)\Rightarrow b&=\inf_{a\in S}(a\Rightarrow b),&
 (\inf S)\Rightarrow b&=\sup_{a\in S}(a\Rightarrow b),\\
 b\Rightarrow(\sup S)&=\sup_{a\in S}(b\Rightarrow a),&
 b\Rightarrow(\inf S)&=\inf_{a\in S}(b\Rightarrow a).
\end{align*}
The two negation equivalences follow similarly from $1-\sup S=\inf_{a\in S}(1-a)$ and its dual. The displayed quantifier-shift rules now yield a prenex equivalent by induction on formulas, after renaming bound variables.
\qed
\end{proof}

Put $\Sigma_0=\Pi_0$ equal to the quantifier-free formulas.  Recursively, $\Sigma_{k+1}$ consists of prenex formulas $\exists\bar x\,\varphi$ with $\varphi\in\Pi_k$, and $\Pi_{k+1}$ consists of prenex formulas $\forall\bar x\,\varphi$ with $\varphi\in\Sigma_k$.  Empty and vacuous blocks are permitted.  Thus,
\begin{equation}
 \Sigma_n\cup\Pi_n\subseteq\Sigma_{n+1}\cap\Pi_{n+1}.
 \label{eq:padding}
\end{equation}
We tacitly use fixed prenex equivalents when assigning a formula to a level of the hierarchy.

\begin{corollary}
\label{cor:hierarchyclosure}
The following closure properties hold over standard semantics, up to logical equivalence.
\begin{enumerate}
\item Negation interchanges $\Sigma_n$ and $\Pi_n$.
\item Each of $\Sigma_n$ and $\Pi_n$ is closed, up to standard equivalence, under finite lattice conjunctions, strong conjunctions, truncated sums, and positive powers.
\item If $\rho\in\Pi_n$ and $q\in\mathbb Q\cap[0,1]$, then $\bar q\to\rho$ is equivalent to a $\Pi_n$-formula and $\rho\to\bar q$ is equivalent to a $\Sigma_n$-formula. Dually, if $\sigma\in\Sigma_n$ and $q\in\mathbb Q\cap[0,1]$, then $\bar q\to\sigma$ is equivalent to a $\Sigma_n$-formula and $\sigma\to\bar q$ is equivalent to a $\Pi_n$-formula.
\end{enumerate}
\end{corollary}

\begin{proof}
Item~(1) follows from Lemma~\ref{lem:prenex}. For item~(2), combine renamed prefixes block by block: the operations $\odot$ and $\oplus$ are continuous and increasing in each
coordinate and hence commute with the relevant suprema and infima.  Item~(3) follows in the same way because $q\Rightarrow(-)$ is increasing, whereas $(-)\Rightarrow q$ is decreasing and therefore reverses every quantifier block. The dual statement must be proved simultaneously in order for the quantifier inductive step to work.
\qed
\end{proof}

For a theory $T$, let $T_{\Pi_n}$ be the set of its standard semantic
$\Pi_n$-consequences.

\begin{definition}
Let $M\subseteq N$.  We call $M$ an \emph{exact substructure} of $N$ if every
atomic $\LangQ(M)$-sentence has the same truth value in $M$ and $N$.

For $n\geq0$,  if $M$ is an exact substructure of $N$ we write $M\preck{n}N$ if for every $\Pi_n$-formula $\rho(\bar x)$
\[
 \val{\rho(\bar a)}{M}=\val{\rho(\bar a)}{N},\qquad \text{for every tuple $\bar a$ from $M$}
\]
\end{definition}

Because negation interchanges the two classes, $M\preck{n}N$ also preserves
the truth degrees of all $\Sigma_n$-formulas with parameters from $M$.  

\begin{definition}
For an $\LangQ$-structure $M$, let $D_n^{\deg}(M)$ contain, for every $\rho(\bar x)\in \Pi_n$, every tuple $\bar a\in M$, and every rational
$q\in\mathbb Q\cap[0,1]$, the sentences
\begin{align*}
 \bar q\to\rho(\bar a)&\quad\text{if }q\leq\val{\rho(\bar a)}{M},\\
 \rho(\bar a)\to\bar q&\quad\text{if }\val{\rho(\bar a)}{M}\leq q.
\end{align*}
\end{definition}

By Corollary~\ref{cor:hierarchyclosure} each displayed formula above is equivalent to a $\Sigma_{n+1}$-formula.  The crisp equality instances in the diagram keep distinct parameter names distinct and ensure that the interpretations of the original constants and functions are preserved.  Thus the map sending $a\in M$ to the interpretation of its name in a model $N$ of the diagram is an embedding; after identifying $M$ with its image,
\begin{equation}
 N\models D_n^{\deg}(M)\quad\Longleftrightarrow\quad M\preck nN,
 \label{eq:diagram}
\end{equation}
where $N$ denotes its reduct to $\LangQ$ on the right-hand side.

\begin{lemma} \label{lem:chain}
Let $n\geq0$, let $\lambda$ be a nonzero ordinal, and let $(M_i)_{i<\lambda}$ be a chain such that $M_i\preck nM_j$ whenever $i<j<\lambda$.  Its exact union $M=\bigcup_{i<\lambda}M_i$ is well defined, and $M_i\preck nM$ for every $i<\lambda$.
\end{lemma}

\begin{proof}
Exactness makes the interpretations of atomic formulas on the union
well defined.  We prove simultaneously, by induction on the remaining
quantifier blocks, that the values of $\Sigma_n$- and $\Pi_n$-formulas with
parameters in $M_i$ agree in $M_i$ and $M$.

The quantifier-free case follows from exactness and induction on the
connectives.  Consider, for example, a formula $\forall x\,\varphi(x,\bar a)$.
By the induction hypothesis, once $b$ belongs to a stage $M_j$, the value of
$\varphi(b,\bar a)$ is the same in all later stages and in the union.  If
$c=\val{\forall x\,\varphi(x,\bar a)}{M_i}$, preservation along the chain gives
\[
 c=\inf_{b\in M_j}\val{\varphi(b,\bar a)}{M_j}
 \quad(j\geq i).
\]
Taking the infimum over all stages yields the same value in $M$.  The
existential case is dual, with suprema in place of infima.
\qed
\end{proof}

\begin{corollary}
\label{cor:varyinglevelchain}
Let $(M_j)_{j<\omega}$ be a chain such that $M_j\preck{j+1}M_{j+1}$ for every
$j<\omega$, and let $M^*=\bigcup_{j<\omega}M_j$ be its exact union.  Then
$M_j\preck{j+1}M^*$ for every $j<\omega$.
\end{corollary}
\begin{proof}
Fix $j<\omega$.  We first show, by induction on $k\geq j$, that
$M_j\preck{j+1}M_k$.  The base case $k=j$ is trivial.  For the inductive
step, suppose $M_j\preck{j+1}M_k$.  By hypothesis
$M_k\preck{k+1}M_{k+1}$, and since $k+1\geq j+1$, repeated use of the
nesting property $\preck{i+1}\Rightarrow\preck i$ gives
$M_k\preck{j+1}M_{k+1}$.  Each relation $\preck{j+1}$ is transitive, being
built from equalities of truth degrees, so
\[
 M_j\preck{j+1}M_k\preck{j+1}M_{k+1}
 \quad\Longrightarrow\quad
 M_j\preck{j+1}M_{k+1}.
\]
Thus $(M_k)_{k\geq j}$ is a $\preck{j+1}$-chain in the sense of
Lemma~\ref{lem:chain}, with exact union
$\bigcup_{k\geq j}M_k=\bigcup_{k<\omega}M_k=M^*$, the two unions agreeing
because the full chain is increasing.  Lemma~\ref{lem:chain} applied to this
$\preck{j+1}$-chain gives $M_j\preck{j+1}M^*$, as required.
\qed
\end{proof}
The following is the degree-sensitive substitute for the classical diagram
argument.

\begin{theorem} \label{thm:extension}
For every $n\geq0$, every theory $T$, and every structure $M$,
\[
 M\models T_{\Pi_{n+1}}
 \quad\Longleftrightarrow\quad
 \text{there is }N\models T\text{ such that }M\preck nN.
\]
\end{theorem}

\begin{proof}
Suppose first that $M\preck nN$ and $N\models T$.  If
$\sigma\in T_{\Pi_{n+1}}$, write
$\sigma=\forall\bar x\,\varphi(\bar x)$ with
$\varphi\in\Sigma_n$.  Since $N\models\sigma$, each instance over $M$ has
value one in $N$, hence also in $M$.  Thus $M\models\sigma$.

Conversely, assume $M\models T_{\Pi_{n+1}}$.  We claim that
$T\cup D_n^{\deg}(M)$ is satisfiable.  Otherwise, Theorem \ref{thm:compactness} gives finitely many members
$\alpha_1(\bar a),\ldots,\alpha_k(\bar a)$ of the degree diagram such that  $T\cup\{\alpha_1(\bar a),\ldots,\alpha_k(\bar a)\}$ is not satisfiable.  
Apply Lemma~\ref{lem:uniformseparation} in the language
expanded by the finitely many parameter names.  There is $m\geq1$ such that
no expansion of a model of $T$ assigns value at least $1-1/m$ to all the
$\alpha_i(\bar a)$.  Replace the parameter names by variables and put
\[
 \delta(\bar x)=
 \bigodot_{i=1}^{k}\alpha_i(\bar x)^{m}.
\]
For every model of $T$ and every tuple $\bar b$, some
$\alpha_i(\bar b)^{m}$ has value zero, by 
Lemma~\ref{lem:power}.  Hence,  $T\models\forall\bar x\,\neg\delta(\bar x)$.
The formula $\delta$ is trivially a
$\Sigma_{n+1}$-formula, so the sentence $\forall\bar x\,\neg\delta(\bar x)$ is a $\Pi_{n+1}$ consequence of $T$. Therefore,  it must have truth value one in $M$.  But every $\alpha_i(\bar a)$ belongs to the degree diagram and has value one in $M$,
so $\delta(\bar a)$ has value one, a contradiction.  A model of
$T\cup D_n^{\deg}(M)$ now gives the required extension by
\eqref{eq:diagram}.
\qed
\end{proof}

For the remainder of this paper, let us put
\[
 \K_n(T)=\Mod(T_{\Pi_{n+1}}).
\]

\begin{corollary} \label{cor:inductive}
The class $\K_n(T)$ is closed under unions of $\preck n$-chains.
\end{corollary}
\begin{proof}
Every axiom in $T_{\Pi_{n+1}}$ has the form
$\forall\bar x\,\varphi$ with $\varphi\in\Sigma_n$.  If it has value one in
every stage, then every instance has value one there.  Lemma~\ref{lem:chain}
preserves the values of these instances in the union, which consequently
satisfies the axiom.
\qed
\end{proof}

\section{Complexity-sensitive infinite forcing}
\label{sec:forcing}
\begin{definition}
Fix $T$ and $n$, and order $\K_n(T)$ by $\preck n$.  For
$M\in\K_n(T)$ and an $\LangQ(M)$-sentence $\varphi$, define its \emph{$n$-forcing
value} recursively by
\begin{align*}
 \forcesval{\varphi}{n}_M&=\val{\varphi}{M}
       &&\text{if }\varphi\text{ is atomic},\\
 \forcesval{\neg\varphi}{n}_M
   &=\inf_{M\preck nN}\bigl(1-\forcesval{\varphi}{n}_N\bigr),\\
 \forcesval{\varphi\to\psi}{n}_M
   &=\inf_{M\preck nN}
       \bigl(\forcesval{\varphi}{n}_N\Rightarrow
             \forcesval{\psi}{n}_M\bigr),\\
 \forcesval{\exists x\,\varphi(x)}{n}_M
   &=\sup_{a\in M}\forcesval{\varphi(a)}{n}_M.
\end{align*}
All extensions in the infima and supremum range in $\K_n(T)$.
\end{definition}
A standard induction gives monotonicity:
\begin{equation}
 M\preck nN\quad\Longrightarrow\quad
 \forcesval{\varphi}{n}_M\leq\forcesval{\varphi}{n}_N.
 \label{eq:monotonicity}
\end{equation}

\begin{definition}
A structure $M\in\K_n(T)$ is \emph{$(T,n)$-generic} if
\[
 \forcesval{\varphi}{n}_M\oplus
 \forcesval{\neg\varphi}{n}_M=1
\]
for every $\LangQ(M)$-sentence $\varphi$.
\end{definition}
Proposition~\ref{prop:generictruth} below shows that $(T,n)$-generic models are exactly those in which the $n$-forcing value and the truth degree coincide.

We now give the generic-extension argument in the relativized order.  Its
key step is the following density statement, adapted from the construction
in \cite[Section~5]{DiNolaGeorgescuSpada2008}.

\begin{lemma}
\label{lem:genericdensity}
Let $A\in\K_n(T)$ and let $\varphi$ be an $\LangQ(A)$-sentence.  There is
$B\in\K_n(T)$ such that $A\preck nB$ and $\forcesval{\varphi}{n}_B\oplus
 \forcesval{\neg\varphi}{n}_B=1$.
\end{lemma}

\begin{proof}
For any $C$ such that $A\preck nC$, let us abbreviate $u_C=\forcesval{\varphi}{n}_C$.  Suppose that the conclusion fails.  Then for every $B\in\K_n(T)$ with $A\preck nB$,
\begin{align*}
 u_B\oplus\forcesval{\neg\varphi}{n}_B
 &=u_B\oplus\inf_{B\preck nC}(1-u_C)\\
 &=\inf_{B\preck nC}(u_C\Rightarrow u_B)<1.
\end{align*}
Here, the second equality follows from the continuity of $u_B\oplus(-)$ on
$[0,1]$.  Hence, some $C\in\K_n(T)$ with $B\preck nC$ satisfies
$u_C\Rightarrow u_B<1$, and therefore $u_B<u_C$.

We construct a $\preck n$-chain
$(B_\xi)_{\xi<\omega_1}$. Set $B_0=A$. At a successor stage, choose
$B_{\xi+1}$ so that
\[
 u_{B_\xi}<u_{B_{\xi+1}}.
\]
At a limit stage, take the exact union. Since every ordinal
$\xi<\omega_1$ is countable, only countable unions occur at limit
stages. By Corollary~\ref{cor:inductive}, every such union belongs to
$\K_n(T)$, and monotonicity gives
\[
 u_{B_\eta}\leq u_{B_\xi}
 \qquad\text{whenever }\eta<\xi<\omega_1.
\]

For each $\xi<\omega_1$, choose, using the density of the rationals, a
rational number $q_\xi$ such that
\[
 u_{B_\xi}<q_\xi<u_{B_{\xi+1}}.
\]
If $\xi<\eta$, then monotonicity gives
$u_{B_{\xi+1}}\leq u_{B_\eta}$, and therefore
\[
 q_\xi<u_{B_{\xi+1}}\leq u_{B_\eta}<q_\eta.
\]
Thus $(q_\xi)_{\xi<\omega_1}$ is a strictly increasing family of
distinct rational numbers, contradicting the countability of
$\mathbb Q$.
\qed
\end{proof}

\begin{theorem}
\label{thm:genericextension}
Every $M\in\K_n(T)$ has an extension $M^*\in\K_n(T)$ such that
$M\preck nM^*$ and $M^*$ is $(T,n)$-generic.
\end{theorem}

\begin{proof}
We begin with a preliminary construction.  Given $A\in\K_n(T)$,
enumerate all $\LangQ(A)$-sentences as $(\varphi_\alpha)_{\alpha<\lambda}$.  Recursively construct a $\preck n$-chain $(A_\alpha)_{\alpha\leq\lambda}$, starting with $A_0=A$.
At a successor stage apply Lemma~\ref{lem:genericdensity} to $A_\alpha$ and $\varphi_\alpha$, and at a limit stage take the exact union.
Corollary~\ref{cor:inductive} ensures that every stage belongs to $\K_n(T)$.  Put $A^+=A_\lambda$.  For each $\alpha<\lambda$, monotonicity shows that the equality obtained at $A_{\alpha+1}$ persists in $A^+$.
Consequently,
\begin{equation}
 \forcesval{\varphi}{n}_{A^+}\oplus
 \forcesval{\neg\varphi}{n}_{A^+}=1
 \quad\text{for every $\LangQ(A)$-sentence $\varphi$.}
 \label{eq:oneround}
\end{equation}

Now, let $M_0=M$ and, for each $j<\omega$, obtain $M_{j+1}=M_j^+$ by one
round of the preceding construction.  Set
\[
 M^*=\bigcup_{j<\omega}M_j.
\]
Again Corollary~\ref{cor:inductive} gives $M^*\in\K_n(T)$ and
$M\preck nM^*$.  Every $\LangQ(M^*)$-sentence contains only finitely many
parameters from $M^*$, so it is already an $\LangQ(M_j)$-sentence for some
$j<\omega$.  Equation~\eqref{eq:oneround} decides it in $M_{j+1}$, and
monotonicity preserves that equality in $M^*$.  Thus $M^*$ is
$(T,n)$-generic.
\qed
\end{proof}

As in the unrelativized forcing construction, genericity is equivalent to the
coincidence of forcing and truth.

\begin{proposition}\label{prop:generictruth}
For $M\in\K_n(T)$, the following are equivalent:
\begin{enumerate}
\item $M$ is $(T,n)$-generic;
\item $\forcesval{\varphi}{n}_M=\val{\varphi}{M}$ for every
      $\LangQ(M)$-sentence $\varphi$.
\end{enumerate}
\end{proposition}

\begin{proof}
The argument is the same as that of \cite[Lemma~5.6]{DiNolaGeorgescuSpada2008}. We provide here a complete proof to show explicitly that only forcing clauses and monotonicity are used.

Assume first that $M$ is generic.  We prove item 2 by induction on $\varphi$. The atomic case is the definition of forcing.  If the assertion holds for $\varphi(x)$, then
\[
 \forcesval{\exists x\,\varphi(x)}{n}_M
 =\sup_{a\in M}\forcesval{\varphi(a)}{n}_M
 =\sup_{a\in M}\val{\varphi(a)}{M}
 =\val{\exists x\,\varphi(x)}{M}.
\]

For negation, notice that $M$ is an extension of itself, so
\[
 \forcesval{\neg\varphi}{n}_M
 \leq 1-\forcesval{\varphi}{n}_M.
\]
Genericity for $\varphi$ gives the reverse inequality, since $\forcesval{\varphi}{n}_M\oplus\forcesval{\neg\varphi}{n}_M=1$. Consequently
\[
 \forcesval{\neg\varphi}{n}_M
 =1-\forcesval{\varphi}{n}_M
 =1-\val{\varphi}{M}.
\]

For implication put
\[
 s=\sup_{M\preck nN}\forcesval{\varphi}{n}_N.
\]
The negation clause and the continuity of $1-x$ give $\forcesval{\neg\varphi}{n}_M=1-s$.  By monotonicity, $\forcesval{\varphi}{n}_M\leq s$, whereas genericity for $\varphi$ gives
\[
 1=\forcesval{\varphi}{n}_M\oplus(1-s),
\]
and hence $s\leq\forcesval{\varphi}{n}_M$. Thus $s=\forcesval{\varphi}{n}_M$. Since $x\Rightarrow b$ is continuous and decreasing in $x$,
\begin{align*}
 \forcesval{\varphi\to\psi}{n}_M
 &=\inf_{M\preck nN}
   \bigl(\forcesval{\varphi}{n}_N\Rightarrow
         \forcesval{\psi}{n}_M\bigr)\\
 &=s\Rightarrow\forcesval{\psi}{n}_M\\
 &=\val{\varphi}{M}\Rightarrow\val{\psi}{M}
  =\val{\varphi\to\psi}{M},
\end{align*}
where the induction hypothesis is used in the last line.  This completes the induction.

Conversely, if forcing and truth coincide for every sentence, then
\[
 \forcesval{\varphi}{n}_M\oplus\forcesval{\neg\varphi}{n}_M
 =\val{\varphi}{M}\oplus(1-\val{\varphi}{M})=1,
\]
so $M$ is generic.  
\qed
\end{proof}

\begin{lemma}
\label{lem:lowagreement}
If $M\in\K_n(T)$ and $\rho$ belongs to $\Sigma_n\cup\Pi_n$, then $\forcesval{\rho}{n}_M=\val{\rho}{M}$.
\end{lemma}
\begin{proof}
Fix a prenex representative of $\rho$.  Every subformula that
occurs in this presentation has complexity at most $n$.
Consequently, whenever $A\preck nB$ and the free variables are assigned
values from $A$, its truth degree is the same in $A$ and $B$.  We prove,
simultaneously for all $A\in\K_n(T)$, that forcing and truth agree on these
subformulas.

The atomic case is immediate.  Suppose first that agreement has been proved
for $\varphi$.  Degree preservation and the induction hypothesis at every
$B$ extending $A$ give
\begin{align*}
 \forcesval{\neg\varphi}{n}_A
 &=\inf_{A\preck nB}\bigl(1-\forcesval{\varphi}{n}_B\bigr)\\
 &=\inf_{A\preck nB}\bigl(1-\val{\varphi}{B}\bigr)
  =1-\val{\varphi}{A}.
\end{align*}
If agreement has been proved for $\varphi$ and $\eta$, then similarly
\begin{align*}
 \forcesval{\varphi\to\eta}{n}_A
 &=\inf_{A\preck nB}
   \bigl(\forcesval{\varphi}{n}_B\Rightarrow
         \forcesval{\eta}{n}_A\bigr)\\
 &=\inf_{A\preck nB}
   \bigl(\val{\varphi}{B}\Rightarrow\val{\eta}{A}\bigr)
  =\val{\varphi}{A}\Rightarrow\val{\eta}{A}.
\end{align*}
Finally,
\[
 \forcesval{\exists x\,\varphi(x)}{n}_A
 =\sup_{a\in A}\forcesval{\varphi(a)}{n}_A
 =\sup_{a\in A}\val{\varphi(a)}{A}.
\]
\qed
\end{proof}

Generic structures are complete with respect to the next existential level.

\begin{proposition}
\label{prop:sigmacomplete}
Let $M$ be $(T,n)$-generic, let $M\preck nN$ with $N\in\K_n(T)$, and let
$\xi$ be a $\Sigma_{n+1}$-sentence of $\LangQ(M)$.  If $N\models\xi$, then
$M\models\xi$.
\end{proposition}

\begin{proof}
Write $\xi=\exists\bar x\,\rho(\bar x)$ with $\rho\in\Pi_n$.  By
Lemma~\ref{lem:lowagreement}, forcing and truth agree on $\rho$ in every
structure in the class, and hence they also agree on the displayed
existential formula.  Since $N\models\xi$,
\[
 \forcesval{\neg\xi}{n}_M
 \leq 1-\forcesval{\xi}{n}_N=0.
\]
Genericity gives $\forcesval{\xi}{n}_M=1$, which is its truth value in $M$.
\qed
\end{proof}

\begin{corollary}
\label{cor:genericconsequences}
Every $(T,n)$-generic structure satisfies $T_{\Pi_{n+2}}$.  In particular, if
$T$ is a $\Pi_{n+2}$-theory, every $(T,n)$-generic structure is a model of
$T$.
\end{corollary}

\begin{proof}
Let $M$ be generic and let
$\sigma=\forall\bar x\,\xi(\bar x)\in T_{\Pi_{n+2}}$, where
$\xi\in\Sigma_{n+1}$.  By Theorem~\ref{thm:extension}, choose
$N\models T$ with $M\preck nN$.  Every instance $\xi(\bar a)$ with
$\bar a\in M$ has value one in $N$, hence in $M$ by
Proposition~\ref{prop:sigmacomplete}.  Taking the infimum over $M$ proves
$M\models\sigma$.
\qed
\end{proof}

\section{The approximate Lablanquie's lemma}
\label{sec:mainlemma}

Lablanquie's lemma \cite[Th\'eor\`eme~2.4 on p.~29]{Lablanquie1978} is a hierarchy-reduction step that connects genericity with the omission of types.  In the notation used here, that result says that if $M$ is a classical
$(T,n)$-generic structure and a $\Sigma_{n+2}$-sentence $\varphi$ of
$\Lang(M)$ is true in $M$, then there is a $\Sigma_{n+1}$-sentence $\psi$ such
that $M\models\psi$ and $T\models\psi\to\varphi$, with no new parameters from $M$: every constant from $M$ occurring in
$\psi$ already occurs in $\varphi$.  
The lemma converts the realization of a formula at level $\Sigma_{n+2}$ into a true antecedent one level lower, which is precisely what produces a lower-level type generating any type realized in a generic model.

The next theorem plays the same role in standard {\L}ukasiewicz semantics, but its conclusion must be adapted to account for truth degrees.  In general, one cannot extract a single exact implication $\psi\to\varphi$.  Instead, for each rational
$r\in[0,1)$ one extracts a true lower-level formula $\psi_r$ entailing the lower bound
$\bar r\to\varphi$.  Since
\[
 \val{\bar r\to\varphi}{N}=1
 \quad\Longleftrightarrow\quad
 r\leq\val{\varphi}{N},
\]
the family of all rational bounds recovers $\val{\varphi}{N}=1$.  

\begin{theorem}[Approximate Lablanquie's lemma]
\label{thm:mainlemma}
Let $M$ be $(T,n)$-generic and let $\varphi$ be a $\Sigma_{n+2}$-sentence of
$\LangQ(M)$ with $M\models\varphi$.  For every rational
$r\in\mathbb Q\cap[0,1)$ there is a
$\Sigma_{n+1}$-sentence $\psi_r$ of $\LangQ(M)$ such that
\begin{enumerate}
\item $M\models\psi_r$;
\item $T\models\psi_r\to(\bar r\to\varphi)$;
\item every constant naming an element of $M$ occurring in $\psi_r$
      already occurs in $\varphi$.
\end{enumerate}
\end{theorem}
\begin{proof}
Let $K$ be the set of all $\Sigma_{n+1}$-sentences of $\LangQ(M)$ having
value one in $M$. The set $K$ has the following maximality property relative
to $T_{\Pi_{n+1}}$: if $\xi\in\Sigma_{n+1}$ and
\begin{equation}\label{eq:maxset}
 K\cup T_{\Pi_{n+1}}\cup\{\xi\}
\end{equation}
is satisfiable, then $\xi\in K$. Indeed, suppose \eqref{eq:maxset} is
satisfiable and let $N$ be a model of it. For every
$\alpha\in D_n^{\deg}(M)$, fix a $\Sigma_{n+1}$-sentence $\alpha^*$
equivalent to $\alpha$. Since $M\models\alpha^*$, we have $\alpha^*\in K$;
therefore $N\models\alpha$. Thus $N\models D_n^{\deg}(M)$, so
$M\preck nN$ by \eqref{eq:diagram}. Moreover, $N\in\K_n(T)$. Since
$N\models\xi$, Proposition~\ref{prop:sigmacomplete} gives $M\models\xi$,
so $\xi\in K$.

Replacing $\varphi$ by a fixed prenex equivalent without introducing new
constants, write $\varphi=\exists\bar x\,\tau(\bar x)$ with
$\tau\in\Pi_{n+1}$. Since the supremum of the values of $\tau$ is one,
choose $\bar a\in M$ such that
\[
 \val{\tau(\bar a)}{M}>r.
\]
Put $\eta=\tau(\bar a)\to\bar r$.  By Corollary~\ref{cor:hierarchyclosure}, this is equivalent to a $\Sigma_{n+1}$-sentence and has value strictly below one in $M$.  By the maximality property, $K\cup T_{\Pi_{n+1}}\cup\{\eta\}$ is not satisfiable.
By Theorem~\ref{thm:compactness} there is a finite $K_0\subseteq K$ such that
\[
 T_{\Pi_{n+1}}\cup K_0\cup\{\eta\}
\]
is not satisfiable.  Apply Lemma~\ref{lem:uniformseparation} to the theory $T_{\Pi_{n+1}}$ and the finite family $K_0\cup\{\eta\}$.  It gives an integer $m\geq1$ such that no model of $T_{\Pi_{n+1}}$ assigns value at least $1-1/m$ to every member of that family.  Put
\[
 \chi=\bigodot_{\kappa\in K_0}\kappa^{m},
\]
with the empty product interpreted as $\top$.  By Corollary~\ref{cor:hierarchyclosure}, this sentence is equivalent to a $\Sigma_{n+1}$-sentence and has value one in $M$.  In any model of
$T_{\Pi_{n+1}}$, either some $\kappa\in K_0$ has value below $1-1/m$, in which case $\kappa^{m}$, and hence $\chi$, has value zero by the formula in the proof of Lemma~\ref{lem:power}, or $\eta$ has value below $1-1/m$, in which case $\neg(\eta^{m})$ has value one.
Consequently,
\[
 T_{\Pi_{n+1}}\models\chi\to\neg(\eta^{m}).
\]
The propositional formula
\[
 \neg\bigl((\alpha\to\bar r)^{m}\bigr)
 \to(\bar r\to\alpha)
\]
is a {\L}ukasiewicz tautology: if $\alpha$ has value at most $r$, its
antecedent has value zero, whereas if $\alpha$ has value greater than $r$,
its consequent has value one. Therefore
\begin{equation}
 T_{\Pi_{n+1}}\models
 \chi\to(\bar r\to\tau(\bar a)).
 \label{eq:bound}
\end{equation}

Write $\chi=\chi(\bar a,\bar b,\bar c)$, where $\bar c$ lists the constants already occurring in $\varphi$ and $\bar b$ lists the remaining parameters. Replace $\bar a,\bar b$ by variables and put
\[
 \psi_r(\bar c)=
 \exists\bar x\exists\bar y\,
 \chi(\bar x,\bar y,\bar c).
\]
This is a $\Sigma_{n+1}$-sentence, and it has value one in $M$.  From \eqref{eq:bound}, monotonicity and continuity of $q\Rightarrow(-)$ with respect to suprema yield
\[
 T\models\psi_r\to(\bar r\to\varphi).
\]
Only the parameters $\bar c$ remain, proving all three assertions.
\qed
\end{proof}

\section{Omitting types theorems}
\label{sec:firstchang}

All formulas in a type $\Gamma(\bar x)$ are assumed to have their free variables among one
fixed finite tuple $\bar x$.

\begin{definition}
A \emph{$T$-type} is a set $\Gamma(\bar x)$ of formulas closed under finite
lattice conjunctions and such that $T\cup \Gamma(\bar x)$ is finitely satisfiable.

We say that a tuple $\bar a\in M$ \emph{realizes} $\Gamma$ in $M$ if $\val{\gamma(\bar a)}{M}=1$ for every $\gamma\in\Gamma$; if no such tuple exists, we say that $M$
\emph{omits} $\Gamma$.  A $T$-type is called a \emph{$\Sigma_j$-type} if each of its members is
equivalent to a $\Sigma_j$-formula.
\end{definition}

\begin{definition}\label{def:existential-type}
Let $\Delta(\bar x)$ and $\Gamma(\bar x)$ be $T$-types.  We say that
$\Delta$ \emph{standard-entails} $\Gamma$, and write
$\Delta\leq_T\Gamma$, if
\[
 T\cup\Delta\models\gamma
 \quad\text{for every }\gamma\in\Gamma.
\]
A $\Sigma_{n+2}$-type $\Gamma$ is \emph{$(n+2)$-existential} if there is no
$\Sigma_{n+1}$-type $\Delta$ such that $\Delta\leq_T\Gamma$.
\end{definition}

\begin{example}
\label{ex:manyvaluedtype}
Let the original language contain a nullary predicate $p$ and a binary
predicate $R$, and put
\[
 T_{1/2}=\{p\to\neg p,\ \neg p\to p\}.
\]
Thus every model of $T_{1/2}$ assigns the non-Boolean value $1/2$ to $p$.
Let $\Gamma_{1/2}(x)$ be the closure under finite lattice conjunctions of the
single formula
\[
 \gamma(x)=p\to\forall y\,R(x,y).
\]
It is a $T_{1/2}$-type: for example, take a one-element structure in which
$p$ and $R$ both have value $1/2$.  Moreover, $\gamma$ is equivalent to the
$\Pi_1$-formula $\forall y\,(p\to R(x,y))$ and hence to a
$\Sigma_2$-formula, too.

We claim that $\Gamma_{1/2}$ is $2$-existential.  Let $\Delta(x)$ be any
$\Sigma_1$-type over $T_{1/2}$.  Theorem~\ref{thm:compactness}, applied after
replacing $x$ by a new constant, gives a model $A\models T_{1/2}$ and
$a\in A$ realizing $\Delta$.  Let $B$ be obtained from $A$ by adding a new element $b$, preserving all atomic values on old tuples, and setting $R(a,b)=0$;
interpret the remaining new atomic instances arbitrarily.  Every
$\Sigma_1$-formula true to degree one of $a$ in $A$ remains true to degree
one in $B$: its old witnesses, or its old family of approximate witnesses,
are still available, and the quantifier-free values on old tuples have not
changed.  The structure $B$ therefore still realizes $\Delta$, whereas
\[
 \val{\gamma(a)}{B}=\tfrac12\Rightarrow0=\tfrac12<1.
\]
Thus $\Delta\not\leq_{T_{1/2}}\Gamma_{1/2}$.  Since $\Delta$ was arbitrary,
$\Gamma_{1/2}$ is $2$-existential.  This example also shows directly how the
first theorem applies outside the Boolean specialization and without adding
truth constants to the original language: a generic exact extension of a
model of $T_{1/2}$ omits $\Gamma_{1/2}$ while $p$ retains value $1/2$.
\end{example}

\begin{proposition}
\label{prop:genericomits}
If $M$ is a $(T,n)$-generic model of $T$, then $M$ omits every
$(n+2)$-existential type.
\end{proposition}

\begin{proof}
Suppose that $\bar a$ realizes an $(n+2)$-existential type $\Gamma$ in $M$.
For every $\gamma\in\Gamma$ and every rational $r\in[0,1)$, apply
Theorem~\ref{thm:mainlemma} to $\gamma(\bar a)$.  We obtain a
$\Sigma_{n+1}$-formula $\psi_{\gamma,r}(\bar x)$ such that
\[
 M\models\psi_{\gamma,r}(\bar a),
 \qquad
 T\models\psi_{\gamma,r}\to(\bar r\to\gamma).
\]
Let $\Delta$ be the closure under finite lattice conjunctions of all these formulas.
The set $\Delta$ is a $\Sigma_{n+1}$-type realized
by $\bar a$.  If $N\models T\cup\Delta$,
then for every $\gamma\in\Gamma$ and every rational $r\in[0,1)$,
\[
 r\leq\val{\gamma}{N}.
\]
It follows that $\val{\gamma}{N}=1$.  Hence
$\Delta\leq_T\Gamma$, a contradiction.
\qed
\end{proof}

\begin{remark}[Comparison with Chang's theorem]
The restriction to $\Sigma_{n+2}$-types is the same quantifier-complexity
restriction that occurs in Chang's classical theorem
\cite{Chang1967}. In Chang's terminology, a
$T$-$(n+2)$-existential type consists of formulas at the
$\Sigma_{n+2}$ level and is not generated by a type at the
$\Sigma_{n+1}$ level. 

There is, however, a difference in the formulation of generation. In
classical logic, compactness of semantic consequence implies that, if a
type $\Delta$ entails a formula $\gamma$, then some finite subset of
$\Delta$ already entails $\gamma$. Since our types are closed under
finite lattice conjunctions, this means that some single $\delta\in\Delta$
satisfies $T\models\delta\to\gamma$.
For standard $[0,1]$-valued semantics, truth-one consequence is not
compact, so in Definition~\ref{def:existential-type} it is necessary to
allow the entire type $\Delta$ to entail $\Gamma$. The rational family
provided by the approximate Lablanquie lemma is the corresponding
many-valued replacement for the single classical implication.

Finally, a single $(T,n)$-generic construction controls formulas only
at the fixed levels $\Sigma_{n+1}$ and $\Sigma_{n+2}$ of the prenex
hierarchy. Types whose formulas have unbounded quantifier complexity
therefore require an iteration through successively higher levels. The
particular class of such types treated by that iteration is introduced
before Theorem~\ref{thm:secondchang}.
\end{remark}

\begin{theorem}[First Chang-style omitting types theorem]
\label{thm:firstchang}
Let $n\geq0$, let $T$ be a $\Pi_{n+2}$-theory, and let $M\models T$.  There is an exact
extension $M^*\models T$ such that
\begin{enumerate}
\item $M^*$ realizes every $\Sigma_{n+1}$-type realized by $M$;
\item $M^*$ omits every $(n+2)$-existential type.
\end{enumerate}
\end{theorem}

\begin{proof}
Since $M\models T$, it belongs to $\K_n(T)$.  By
Theorem~\ref{thm:genericextension}, choose a $(T,n)$-generic $M^*$ with
$M\preck nM^*$.  Corollary~\ref{cor:genericconsequences} gives
$M^*\models T$, and Proposition~\ref{prop:genericomits} gives item 2.

For item 1, let $\bar a\in M$ realize a $\Sigma_{n+1}$-type $\Gamma$, and
fix $\gamma(\bar x)\in\Gamma$. Choose for $\gamma(\bar a)$ an equivalent
representative $\exists\bar y\,\rho(\bar a,\bar y)$ with
$\rho\in\Pi_n$. Its truth value in $M$ is the supremum of the values of the
instances $\rho(\bar a,\bar b)$ as $\bar b$ ranges over $M$. These values are
preserved in $M^*$, and all the old approximate witnesses remain available
there. Hence the corresponding supremum in $M^*$ is at least the supremum in
$M$, and therefore equals one. Thus $M^*\models\gamma(\bar a)$. Since
$\gamma\in\Gamma$ was arbitrary, $\bar a$ realizes $\Gamma$ in $M^*$.
\qed
\end{proof}

The second theorem uses genericity at successively higher levels.  We first
isolate the semantic non-generation condition needed in its statement.  If
$\Delta(\bar x)$ is a $T$-type and $m<\omega$, put
\[
 \Delta_{\Sigma_m}=
 \{\delta\in\Delta:\text{$\delta$ is equivalent over $T$ to a
 $\Sigma_m$-formula}\}.
\]

\begin{definition}\label{def:infty-universal}
A $T$-type $\Gamma(\bar x)$ is \emph{infinity-universal} if it has an
enumeration $(\gamma_k)_{k<\omega}$ and a strictly increasing function
$f\colon\omega\to\omega$ such that
\begin{enumerate}
\item $T\models\gamma_{k+1}\to\gamma_k$ for every $k$;
\item $\gamma_k\in\Pi_{f(k)}$ and $\gamma_k$ is not equivalent over $T$
      to any $\Sigma_{f(k)}$-formula, for every $k$;
\item there do not exist a $T$-type $\Delta(\bar x)$ and $k_0<\omega$ such
      that
      \[
       T\cup\Delta_{\Sigma_{f(k)}}\models\gamma_k
       \qquad\text{for every }k\geq k_0.
      \]
\end{enumerate}
\end{definition}

Thus clause~(1) gives a descending basis for the type, clause~(2) requires
genuinely unbounded quantifier complexity, and clause~(3) says that the tail
of this basis is not generated by the corresponding lower-complexity
fragments of a single type.  In particular, the definition is independent
of the forcing construction and of the particular syntactic form of the
approximate Lablanquie lemma.

\begin{remark}[Relation with the classical definition]
In the Boolean specialization, Definition~\ref{def:infty-universal} is
equivalent to Lablanquie's classical condition
\cite[D\'efinition~3.3 and Th\'eor\`eme~3.4]{Lablanquie1978}.  Indeed, compactness of classical semantic consequence implies that
\[
 T\cup\Delta_{\Sigma_{f(k)}}\models\gamma_k
\]
holds if and only if some finite subset of
$\Delta_{\Sigma_{f(k)}}$ entails $\gamma_k$ over $T$.  Since types are closed
under finite conjunctions and $\Sigma_{f(k)}$ is closed under such
conjunctions, this is equivalent to the existence of
$\delta_k\in\Delta$, equivalent over $T$ to a $\Sigma_{f(k)}$-formula, such
that
\[
 T\models\delta_k\to\gamma_k.
\]

For standard $[0,1]$-valued semantics, truth-one consequence is not compact,
so the entire fragment $\Delta_{\Sigma_{f(k)}}$ cannot in general be
replaced by a single formula.  The approximate Lablanquie lemma nevertheless
shows that, if $\gamma_k$ is realized in the iterated generic extension,
this fragment entails $\gamma_k$: it supplies, for every rational $r<1$, a
true formula of complexity $\Sigma_{f(k)}$ forcing the value of $\gamma_k$
to be at least $r$.  
\end{remark}

\begin{theorem}[Second Chang-style omitting types theorem]
\label{thm:secondchang}
For every theory $T$ and every $M\models T$, there is an exact extension
$M^*\models T$ omitting every infinity-universal $T$-type.
\end{theorem}

\begin{proof}
Construct a chain $(M_j)_{j<\omega}$ as follows.  Let $M_0$ be a $(T,0)$-generic extension of $M$.  Having constructed $M_j$, note from Corollary~\ref{cor:genericconsequences} that $M_j\models T_{\Pi_{j+2}}$, so $M_j\in\K_{j+1}(T)$.  Choose a
$(T,j+1)$-generic $M_{j+1}$ with
\[
 M_j\preck{j+1}M_{j+1}.
\]
Let $M^*=\bigcup_{j<\omega}M_j$.

We first verify that $M^*\models T$.  Given $\sigma\in T$, choose $j$ so large that $\sigma$ is equivalent to a $\Pi_{j+1}$-sentence.  The structure $M_j$ satisfies $\sigma$ by Corollary~\ref{cor:genericconsequences}.  The chain $(M_i)_{i<\omega}$ satisfies the hypotheses of Corollary~\ref{cor:varyinglevelchain}, so $M_j\preck{j+1}M^*$; in particular, this relation preserves the value of $\sigma$ in $M^*$.

Suppose now that a tuple $\bar a\in M^*$ realizes an infinity-universal type
$\Gamma=(\gamma_k)_{k<\omega}$, witnessed by a strictly increasing function
$f$.  Let
\[
 \Delta(\bar x)=
 \{\delta(\bar x):M^*\models\delta(\bar a)\}
\]
be the truth-one type of $\bar a$ in $M^*$.  Choose $s$ with
$\bar a\in M_s$.  For all sufficiently large $k$, put
$j=f(k)-1\geq s$. Then $M_j$ is $(T,j)$-generic.
Corollary~\ref{cor:varyinglevelchain} gives $M_j\preck{j+1}M^*$, and since
$f(k)=j+1$, this is $M_j\preck{f(k)}M^*$.  Since
$\gamma_k\in\Pi_{f(k)}$ and $M^*\models\gamma_k(\bar a)$, the definition of
$\preck{f(k)}$ gives directly $M_j\models\gamma_k(\bar a)$.

Regard $\gamma_k$ as a $\Sigma_{f(k)+1}$-formula and apply the approximate Lablanquie's lemma at level $j=f(k)-1$.  For each rational $r\in[0,1)$ it supplies a $\Sigma_{f(k)}$-formula $\delta_{k,r}(\bar x)$ such that
\[
 M_j\models\delta_{k,r}(\bar a),
 \qquad
 T\models\delta_{k,r}\to(\bar r\to\gamma_k).
\]
By the parameter clause in Theorem~\ref{thm:mainlemma}, after replacing the
names of $\bar a$ by the variables $\bar x$, the formula
$\delta_{k,r}(\bar x)$ has no additional parameters.  Its degree at
$\bar a$ is preserved in $M^*$, so $\delta_{k,r}\in\Delta_{\Sigma_{f(k)}}$.
Now let $N\models T$ and suppose that a tuple $\bar b\in N$ realizes
$\Delta_{\Sigma_{f(k)}}$.  For every rational $r<1$ we then have
\[
 N\models\delta_{k,r}(\bar b)
 \quad\text{and hence}\quad
 \val{\gamma_k(\bar b)}{N}\geq r.
\]
It follows that $N\models\gamma_k(\bar b)$.  Consequently,
\[
 T\cup\Delta_{\Sigma_{f(k)}}\models\gamma_k
\]
for every sufficiently large $k$.  Thus $\Delta$ eventually generates the
tail of $\Gamma$, contradicting clause~(3) of
Definition~\ref{def:infty-universal}.  Hence $M^*$ omits $\Gamma$.
\qed
\end{proof}
\section*{Acknowledgements}
The author expresses his gratitude to Antonio Di Nola and George Georgescu for useful discussions on forcing and its connections with omitting types theorems.

\end{document}